# Integral Regular Truncated Pyramids with Rectangular Bases

*Konstantine Zelator*

*Department of Mathematics*

*301 Thackeray Hall*

*University of Pittsburgh*

*Pittsburgh, PA 15260, U.S.A.*

*Also: Konstantine Zelator*

*P.O. Box 4280*

*Pittsburgh, PA 15203, U.S.A.*

*Email addresses: 1) konstantine_zelator@yahoo.com*
*2) spaceman@pitt.edu*

## 1. Introduction

A regular truncated pyramid with rectangular bases can be obtained by intersecting a regular pyramid with a rectangular base; with a plane parallel to the base of the pyramid, somewhere between the base and the apex of the pyramid. Consequently, a regular truncated pyramid with rectangular bases; consists of two parallel rectangular bases such that the (straight) line connecting the two centers (of the two rectangular bases) is perpendicular to the two parallel planes containing the two rectangular bases; that line being one of the axes of symmetry of the (truncated) pyramid. The other four faces of the truncated pyramid, are all isosceles trapezoids; occurring in two pairs, each pair containing two congruent isosceles trapezoids facing (or opposite) each other.

**Definition 1:** *An integral regular truncated pyramid with rectangular bases is a regular truncated pyramid with rectangular bases whose twelve edges have integer lengths; whose height (the distance between the two bases) is also an integer; and whose volume is an integer as well.*

Of course in the case of a regular truncated pyramid with rectangular bases; the four lateral edges have the same length $t$. The bottom (larger) base is a rectangle of dimensions $a$ and $b$; and the top base a rectangle of dimensions $c$ and $d$. So if $H$ is the height of the pyramid; we really have six lengths involved: $t, H, a, b, c,$ and $d$. Because we are dealing with a regular pyramid, the two rectangle centers are aligned perpendicularly (to the two parallel bases) ; and thus, the two rectangular bases are either both nonsquare rectangles; which means that (with larger $a$ of $a$ and $b$; $c$ the larger of $c$ and $d$) we have, $a>b$, $c>d$, $a>c$, and $b>d$. Or, alternatively, both bases are squares in which case we *have* $a=b > c=d$.

The aim of this paper is to parametrically describe the entire set of integral regular truncated pyramids with rectangular bases. We do so in **Proposition 1, Section 7**.

In **Section 2** we provide some illustrations of a regular truncated pyramid's (with rectangular bases) faces and cross sections. In **Section 3**, we derive a formula for the volume of such a pyramid. In **Section 4**, we derive the key equation of this paper. It is this equation that allows us to parametrically determine the set of all such integral pyramids. The key equation reduces to the diophantine equation, $t^2=x^2+y^2+z^2$.

This is a well-known equation, and it's general solution (in positive integers $x$, $y$, $z$, and $t$) can be found in W. Sierpinski's book *Elementary Theory of Numbers.* (see **Reference [1]**)

Also note that since $x, y, z, t$ are positive integers in the above equation; and $\sqrt{3}$ is an irrational number; at most two among $x,y,z$ can be equal. Accordingly, we distinguish between the case when $x, y, z$ are distinct; and the case when two among them are equal. The latter case leads to the diophantine equation $Z^2=X^2+2Y^2$, whose general solution we state in **Section 6**. In **Section 5**, we state the general positive integer solution of the equation $t^2=x^2+y^2+z^2$.

In **Section 8, Proposition 2**, we describe a particular 3-paremeter set or family of integral regular truncated pyramids with rectangular bases. In **Section 9, Proposition 3**, we give a 4-parameter description of the set of integral regular truncated pyramids with square bases. We end this paper with some closing remarks in **Section 10**.

## 2. Illustrations

In the illustrations below, $a$ and $b$ are the sidelengths of the bottom rectangular base; with $a \geq b$. And $c$ and $d$ the sidelengths of the top rectangular base; with $c \geq d$. And consequently, $a > c$ and $b > d$ must hold as well. And $t$ stands for the common length of the four lateral edges.

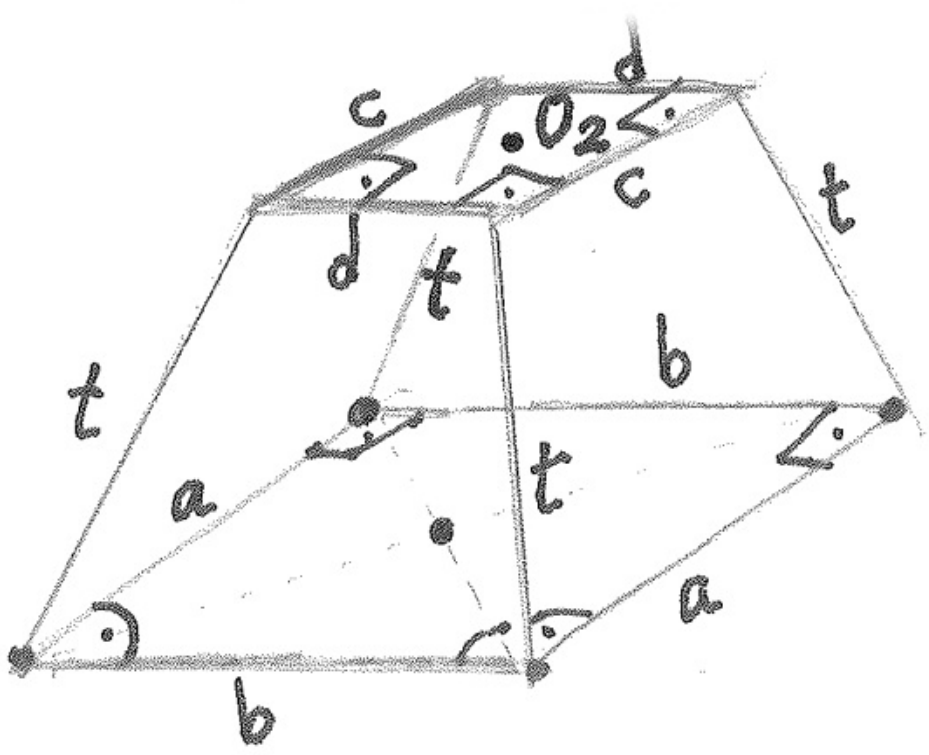

$O_1$ is the center of the bottom rectangular base

$O_2$ the center of the bottom base; the straight line segment $\overline{O_1O_2}$ is perpendicular (or orthogonal) to the two parallel bases.

**Figure 1: 3-D picture of a truncated regular pyramid with rectangular bases**

**Figure 2: Perpendicular (or orthogonal) projection of the top base onto the bottom page**

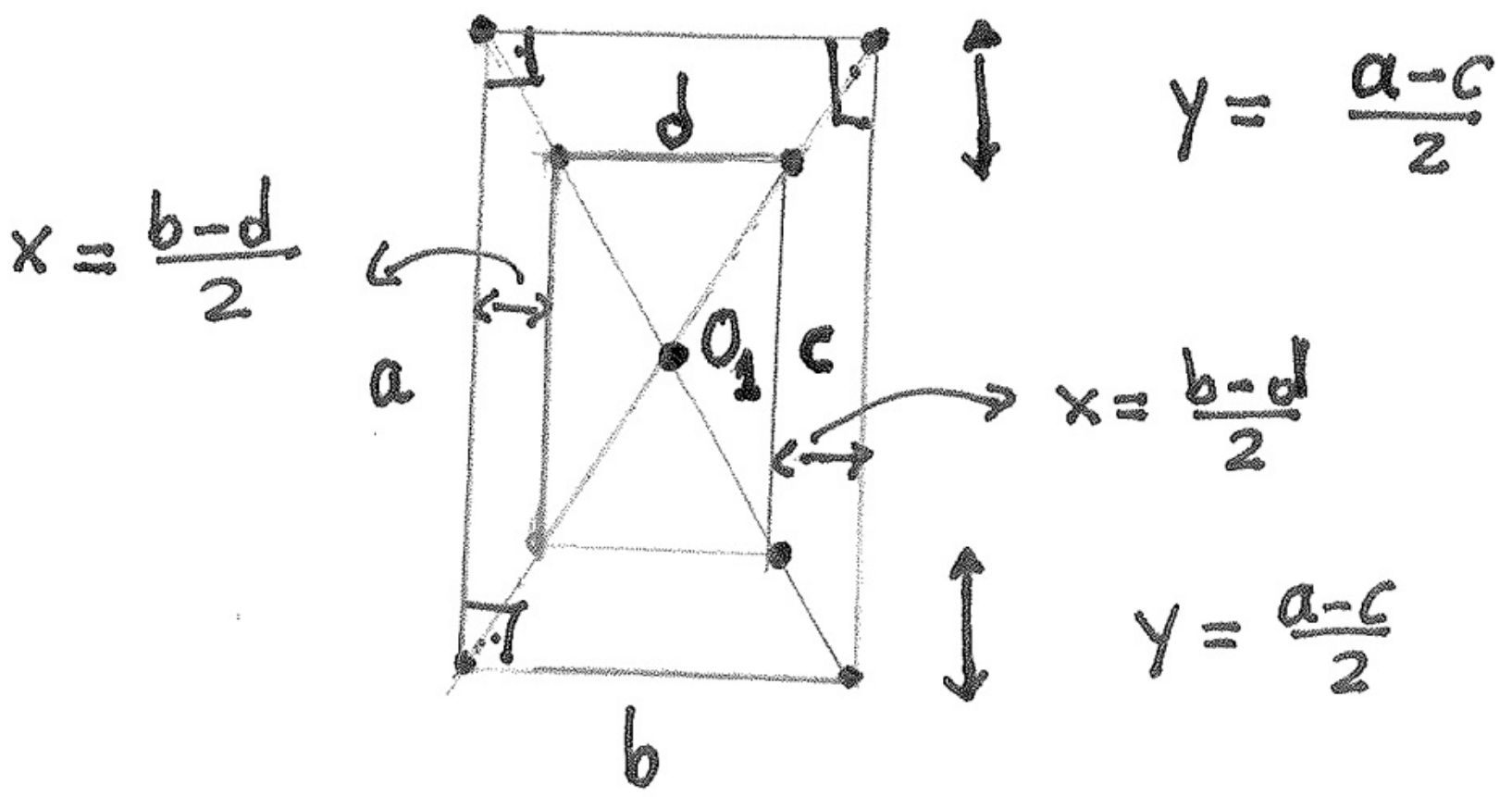

**Figure 3: Bottom base**

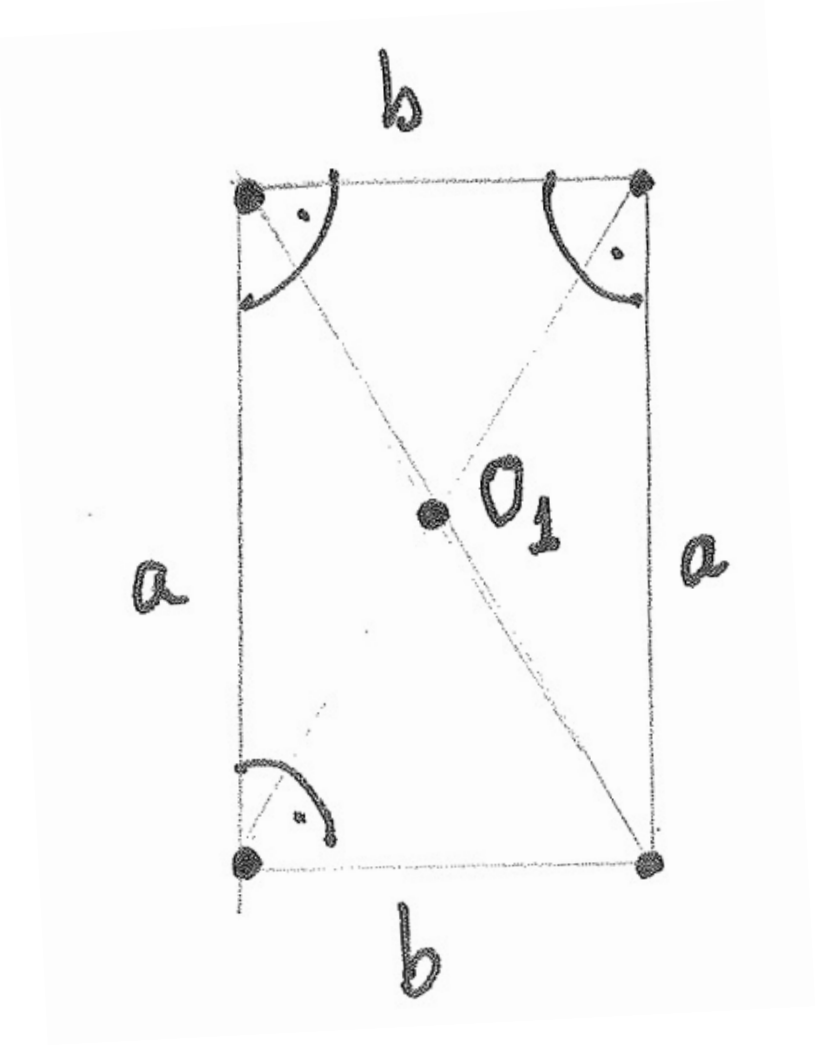

**Figure 4: Top base**

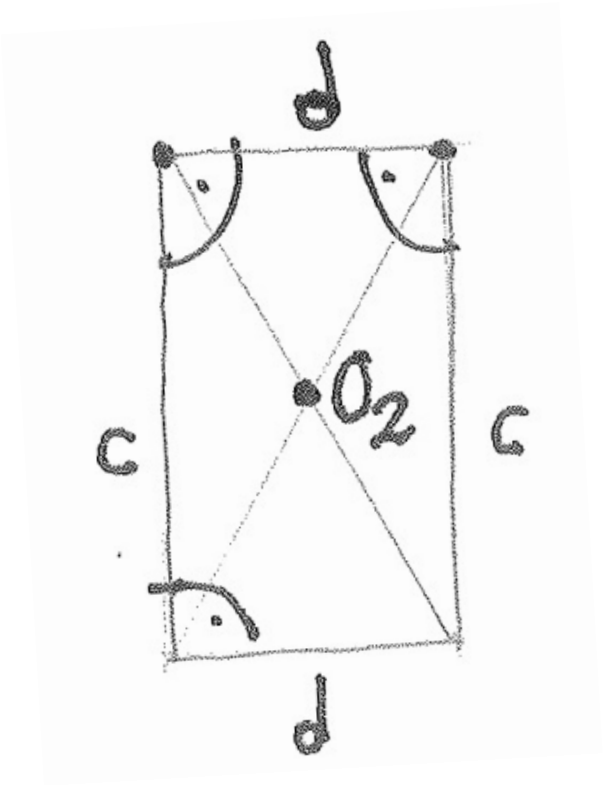

**Figure 5: Two congruent lateral faces (isosceles trapezoids)**

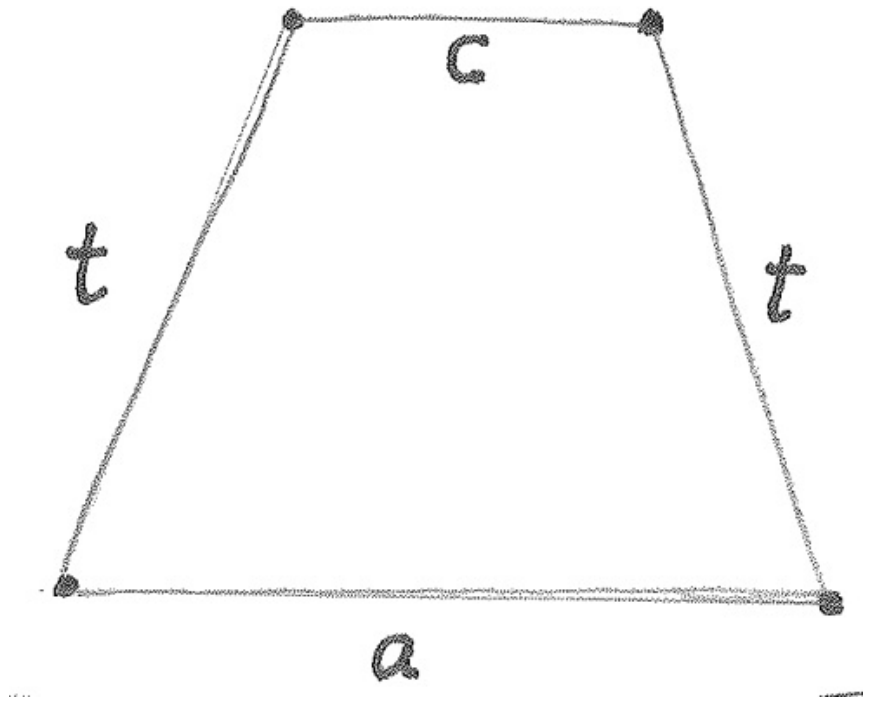

**Figure 6: Two congruent lateral faces (isosceles trapezoids)**

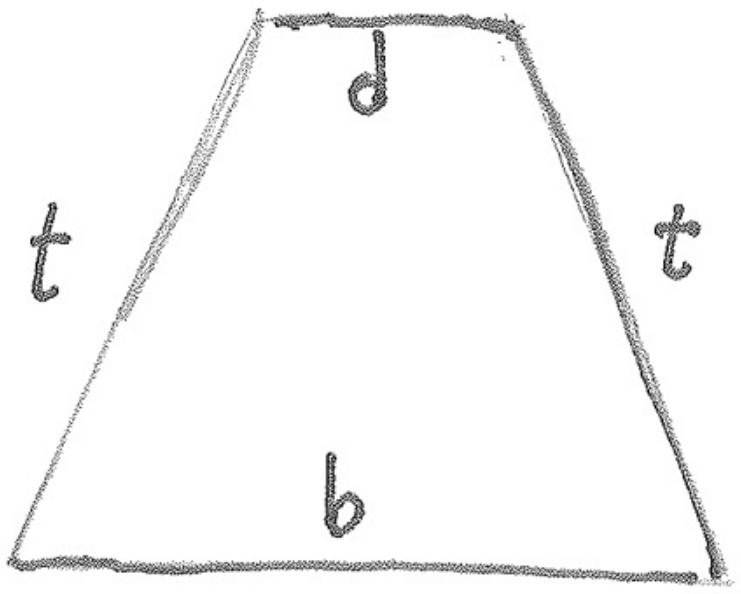

3. **Geometric Considerations**

**Figure 7***: a≥b, c≥d, a>c, b>d*

$$\left|\overline{P_1P_3}\right| = b, \left|\overline{Q_1Q_3}\right| = d,$$
$$\left|\overline{P_1P_2}\right| = a, \left|\overline{Q_1Q_2}\right| = c$$

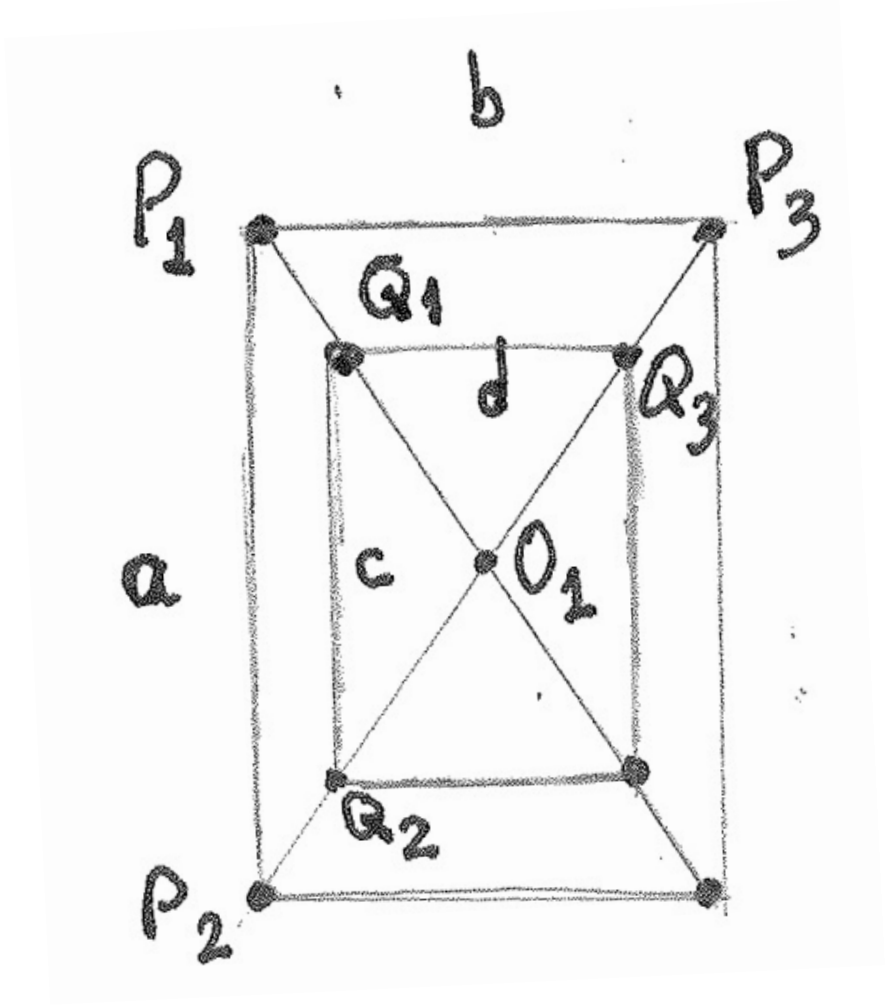

A. **The Key Ratio**

From the similar triangle $O_1Q_1Q_3$ and $O_1P_1P_3$ we have,

$$\frac{\left|\overline{O_1Q_1}\right|}{\left|\overline{O_1P_1}\right|} = \frac{\left|\overline{Q_1Q_3}\right|}{\left|\overline{P_1P_3}\right|} = \frac{d}{b}$$

And from the similar triangles $O_1Q_1Q_2$ and $O_1P_1P_2$ we obtain,

$$\frac{\left|\overline{O_1Q_1}\right|}{\left|\overline{O_1P_1}\right|} = \frac{\left|\overline{Q_1Q_2}\right|}{\left|\overline{P_1P_2}\right|} = \frac{c}{a}$$

Therefor we conclude that $\frac{d}{b} = \frac{c}{a}$ or equivalently,

$$\frac{a}{b} = \frac{c}{a} \qquad \textbf{(1)}$$

This being the key ratio; an important condition that the four base sidelengths of a regular truncated pyramid with rectangular bases must satisfy.

Now, consider the case in which *a,b,c,d* are positive integers. In number theory, Euclid's Lemma postulates that if a positive integer $i_1$ divides the product $i_2 i_3$ of two other positive integers $i_2$ *and* $i_3$ *;* and $i_1$ is relatively prime to $i_2$. Then $i_1$ must be a divisor of $i_3$*.* Using Euclid's Lemma in **(1)** (one rewrites **(1)** in the form *ad=bc*) one can easily prove that **(1)** implies,

$$\left\{\begin{array}{l} a = Nk_1,\ b = Nk_2,\ c = Mk_1,\ d = Mk_2 \\ \textit{where } N,\ M,\ k_1,\ k_2 \textit{ are positive integers} \\ \textit{with } k_1,\ k_2 \textit{ being relatively prime}; \\ gcd\ (k_1,\ k_2)\ =\ 1 \end{array}\right\} \quad \textbf{(2)}$$

The conditions in **(2)** can be understood as lowest terms conditions. Note that *N* is none other than the greatest common of *a* and *b*; *N=gcd (a,b).* Likewise *M=gcd (c,d).*

Observe that since *a>c* and *b>d; a* and *b* cannot be relatively prime (for this would imply *N=1*; and consequently *a≤c* and *b≤d*; contrary to *a>c* and *b>d).*

*On the other hand b and d can be relatively prime:*

*When c and d are relatively prime in* ***(2)****; we have M=1 and so* $c=k_1$*,* $d=k_2$*. Thus in such a case, a=Nc and b=Nd.*

*There is also the square case: a=b and c=d. When a=b and c=d; it follows that* $k_1 = k_2$*; which in turn implies (since* $gcd\ (k_1,k_2)=1$*) that* $k_1=k_2=1$*. So in the square case a=b=N and c=d=M ; with the key ratio being equal to 1.*

### B. A Formula for the Volume of a Regular Truncated Pyramid with Rectangular Bases

The straight line $\overleftrightarrow{O_1O_2}$ that connects the centers $O_1$ and $O_2$ of the two rectangular bases is perpendicular or orthogonal to the two parallel planes that contain the two rectangular bases. Consequently, the line $\overleftrightarrow{O_1O_2}$[; as well the four lines containing the four lateral edges; they all meet at an apex point *A,* as illustrated in **Figure 8** below.

3-D picture

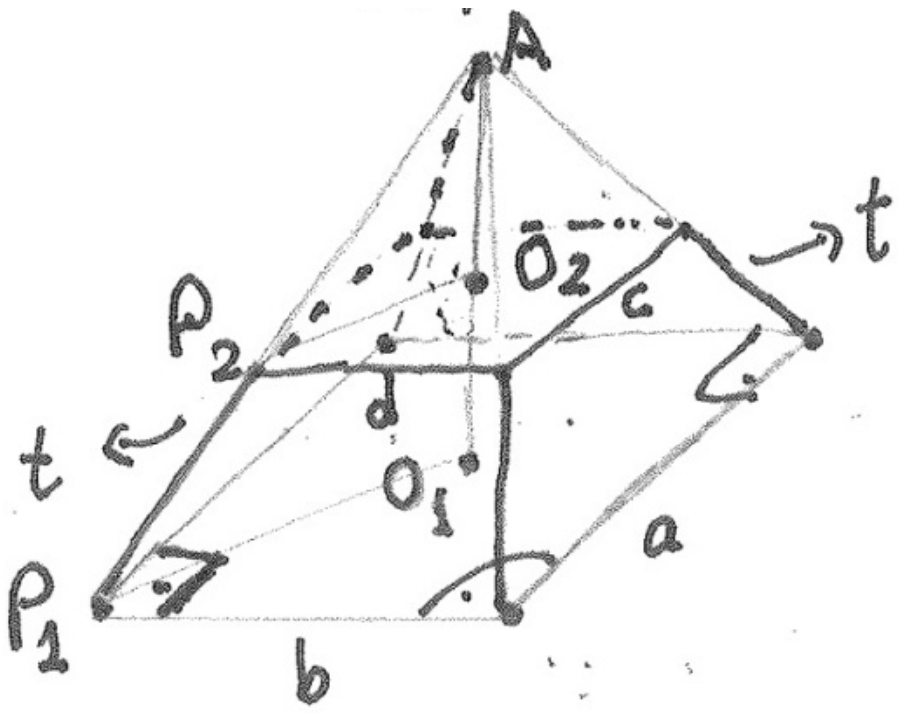

**Figure 8**

**$H_2$=height of the top pyramid**

**$H_1$=height of the larger pyramid**

**$H=H_1-H_2=\left|\overline{O_1O_2}\right|$ = height of the truncated pyramid**

Let *V* be the volume of the truncated pyramid, $V_1$ the volume of the larger pyramid; and $V_2$ the volume of the larger pyramid; and the $V_2$ the volume of the top or smaller pyramid. Then,

$$\left\{V=V_1-V_2=\frac{1}{3}H_1ab-\frac{1}{3}H_2cd\right\} \quad \textbf{(3)}$$

From the similar right (or 90 degree) triangles $P_2O_2A$ and $P_1O_1A$ it is clear that,

$$\frac{H_1}{H_2}=\frac{\left|\overline{P_1O_1}\right|}{\left|\overline{p_2O_2}\right|}=\frac{\left(\sqrt{a^2+b}\,{}^{2}\!\big/2\right)}{\left(\sqrt{c^2+d}\,{}^{2}\!\big/2\right)}=\frac{\sqrt{a^2+b^2}}{\sqrt{c^2+d^2}};$$

$$\frac{H_1}{H_2}=\frac{b\sqrt{\left(\frac{a}{b}\right)^2+1}}{\mathrm{d}\sqrt{\left(c\!\big/d\right)^2+1}}$$ ; and by **(1)** it follows that

$$\frac{H_1}{H_2}=\frac{b}{d} \quad \textbf{(4)}$$

From **(4)** we further get,

$$\frac{H_1 - H_2}{H_2} = \frac{b-d}{d} = \frac{a-c}{c}; \qquad \textbf{(5)}$$

**by (1)**

And using $H = H_1 - H_2$ we further get,

$$H_2 = \frac{Hc}{a-c} \qquad \textbf{(6)}$$

And by **(6)** and **(4)**, $H_1 = \frac{Ha}{a-c}$ **(7)**

Combining **(7), (6)** and **(3);** yields

$V = \frac{1}{3} \cdot \frac{H}{a-c} \cdot \left[ a^2 b - c^2 d \right]$; and using b=$\frac{ad}{c}$; further gives

$V = \frac{H \cdot \left(a^3 - c^3\right) \cdot d}{3c\left(a-c\right)}$; which combined with

$a^3 - c^3 = \left(a-c\right)\left(a^2 + ac + c^2\right)$ ; produces

$$\left\{ V = \frac{d \bullet H \bullet (a^2 + ac + c^{2)}}{3c} \right\} \qquad \textbf{(8)}$$

Now, when *d,H,a*, and *c* are positive integers; the volume *V* is also an integer if and only if *3c* is a divisor of the product $d \cdot H \cdot \left(a^2 + ac + c^2\right)$ Observe that if *a=c(mod3)*. Then,

$$a^2 + ac + c^2 \equiv a^2 + a^2 + a^2 \equiv 3a^2 \equiv O\left(\mathrm{mod}\,3\right) .$$

So, when *a* and *c* belong to the same congruence class modulo 3;

The integer $a^2 + ac + c^2$ is divisible by 3.

## 4. The Key Equation

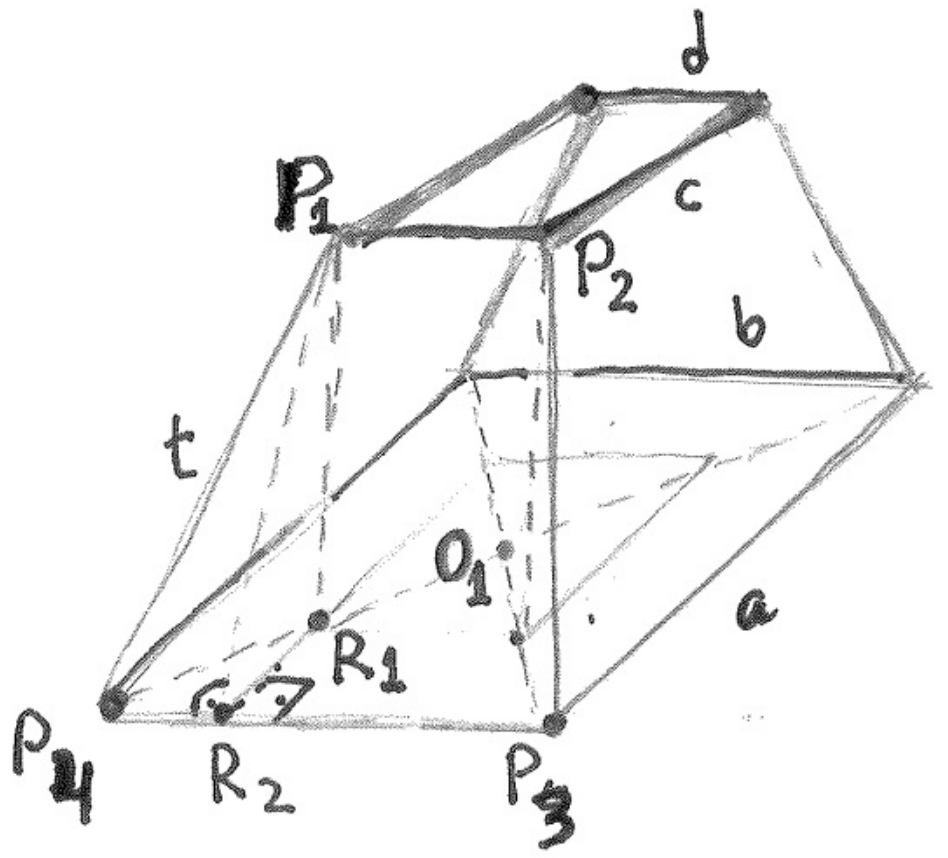

**Figure 9:** $\overline{P_1P_2}$ **is one (of the four) edges of the top rectangular base;** $\overline{P_3P_4}$ **is one (of the four) edges of the bottom rectangular base;** $P_1P_2P_3P_4$ **is an isosceles trapezoid; one of the four lateral faces. The point** $R_1$ **is the orthogonal projection of** $P_1$ **onto the bottom base;** $R_2$ **is the perpendicular projection of** $R_1$ **onto the edge** $\overline{P_3P_4}$ **; and also of** $P_1$ **onto the edge** $\overline{P_3P_4}$

As we have already seen, $\left|\overline{P_1P_4}\right| = t$ , $\left|\overline{P_1R_1}\right| = H =$ height of the truncated pyramid; $\left|\overline{R_1R_2}\right| = y = \dfrac{a-c}{2}$ and $\left|\overline{P_4R_2}\right| = x = \dfrac{b-d}{2}$ (**see Figure 2**). Thus, from the right triangles $P_1R_2P_4$ and $P_1R_1R_2$ we have;

$$\left\{\begin{array}{l} \left|\overline{P_1P_4}\right|^2 = \left|\overline{P_4R_2}\right|^2 + \left|\overline{P_1R_2}\right|^2 \\ and \left|\overline{P_1R_2}\right|^2 = \left|\overline{R_1R_2}\right|^2 + \left|\overline{P_1R_1}\right|^2 \end{array}\right\} \qquad \textbf{(9)}$$

And so from **(9)** we obtain the key equation

$$t^2 = H^2 + \left(\frac{a-c}{2}\right)^2 + \left(\frac{b-d}{2}\right)^2$$ ; or equivalently.

$$\left\{4t^2 = 4H^2 + (a-c)^2 + (b-d)^2\right\} \qquad \textbf{(10)}$$

So, when *t, H, a, b, c, d* are positive integers; with *a>c, b>d, a≥b, and c≥d*. The sum of the squares $(a-c)^2 + (b-d)^2$ must be according to **(10)** ; congruent to *0*(mod4). But the square of an integer is congruent to *0* or *1* modulo 4; according to whether that integer is even or odd. It becomes clear that both positive integers, *a-c* and *b-d* must be even. Therefore, in conjunction with **(10)** we must have,

$$\left\{\begin{array}{l} t^2 = H^2 + y^2 + x^2 \\ a-c = 2y,\ b-d = 2x; \\ with\, t, H, x, y, a, c, b, d \\ all\ being\ \ positive\ integers \\ that\ also\ satisfy \\ a \geq b\ and\ c \geq d \end{array}\right\}$$

Note that the conditions *a>c* and *b>d* are automatically satisfied in view of the fact that *y* and *x* are positive integers. Also *a≥b* $\Leftrightarrow\ c + 2y \geq d + 2x$; $c - d \geq 2(x-y)$ But also from **(1) ,** we have

$\frac{a}{c} = \frac{b}{d} \Leftrightarrow \frac{a-c}{c} = \frac{b-d}{d}$ ; and so, $\frac{2y}{c} = \frac{2x}{d}$; $c = \frac{y \cdot d}{x}$ And so, $c - d = \frac{yd}{x} - d = \frac{\mathrm{d}(y-x)}{x}$

And so, the conditions $\left\{\begin{array}{l} c - d = \frac{\mathrm{d}(y-x)}{x} \\ And\ \ c - d \geq 2(x-y) \end{array}\right\}$; further require,

$$\frac{\mathrm{d}(y-x)}{x} \geq 2(x-y);$$

$$\frac{\mathrm{d}(y-x)}{x} + 2(y-x) \geq 0;$$

$$\frac{(y-x)}{x} \cdot [d + 2x] \geq 0;$$

which is equivalent (since *d, x, y* are positive integers) to

$y - x \geq 0$; or $x \leq y$

We can finalize our conditions as follows.

$$\left\{\begin{array}{l} t^2 = H^2 + y^2 + x^2 \\ a - c = 2y, b - d = 2x, and, y \geq x; \\ Where\ t,\ H,\ y,\ x,\ c,\ d\ are\ positive \\ integers;\ and\ thus\ so\ are\ a\ and\ b. \\ (And\ thus\ also,\ the\ conditions\ a > c,\ b > d,\ a \geq b,\ c \geq d) \\ \\ And\ with\ \mathrm{c} = \frac{y \cdot d}{x} \left(so\ x\ being\ a\ divisor\ \mathrm{of}\ y \cdot d\ \right) \\ \left(And\ thus\ the\ necessary\ condition\ (1),\ \frac{a}{c} = \frac{b}{d} is\ also\ satisfied\right) \end{array}\right\} \quad \mathbf{(11)}$$

**5. The Diophantine Equation** $t^2 = 2^2 + y^2 + x^2$

A derivation and the general solution in positive integers, to the equation $t^2 = z^2 + y^2 + x^2$ **;** can be found in W. Sierpinski's book *Elementary Theory of Numbers* (see **Reference [1]**). The entire solution set in positive integers *x, y, z, t*; can be described in terms of three integer parameters *m, l,* and *n* as follows.

$$\left\{\begin{array}{l} The\ positive\ integers\ x,\ y,\ z,\ and\ t\ satisfy\ the\ equation \\ t^2 = z^2 + y^2 + x^2, if\ and\ only\ if \\ t = \dfrac{l^2+m^2+n^2}{n} = \dfrac{l^2+m^2}{n} + n, z = 2m, y = 2l, x = \dfrac{l^2+m^2-n^2}{n} \\ = \dfrac{l^2+m^2}{n} - n \\ Where\ m,\ n,\ l\ are\ positive\ integers\ such\ that \\ 1 \le n < \sqrt{l^2+m^2}\ \left(so\ that\ n^2 < l^2+m^2\right) And\ with\ n\ being \\ a\ divisor\ of\ l^2+m^2 \end{array}\right\} \quad \textbf{(12)}$$

In **Reference [1],** the reader will also find a list of particular solutions to the above equation.

Observe that when $l^2+m^2$ is a prime number (as it is well known in number theory; only 2 and primes congruent to 1 modulo 4, can be expressed as sums of integer squares. Primes congruent to 3 mod 4 cannot) ; n=1 is the only positive divisor of $l^2+m^2$, which is less than $\sqrt{l^2+m^2}$ ; since the only positive divisors of a prime; are itself and 1. Regardless of whether $l^2+m^2$ is a prime or not, the choice *n=1* generates a subfamily of solutions. We have the following.

*A subfamily of the entire family of solutions of positive integer solutions to the equation in* ***(12)*** *is the set,* $s = \left\{\left(x,y,z,t\right) \middle| x = l^2+m^2-1,\ y = 2l,\ \ z = 2m,\ t = l^2+m^2+1,\ l,\ m \in\ \mathbb{Z}^+\right\}$ **(13)**

So, the set *S* is a 2-parameter subfamily of positive integer solutions to the diophantine equation in **(12).** Observe that at most two of x,y,z can be equal; since $\sqrt{3}$ is an irrational number. When two of *x, y, z* to are equal. The equation $t^2 = z^2 + y^2 + x^2$ reduces to the Diophantine equation $Z^2 = X^2 + 2Y^2$ .

## 6. The Diophantine equation $Z^2 = X^2 + 2Y^2$

This equation is well-known in the literature. For this and related historical material, see **Reference [2].** In **Reference [3]** (this is an article published in 2006), the reader can find a detailed analysis and derivation of the general positive integer solution to the equation $Z^2 = X^2 + k \cdot Y^2$, where *k* is a fixed positive integer, $k \geq 2$. We state the general solution to $Z^2 = X^2 + 2Y^2$ below.

All the solutions, in positive integers *X, Y, Z* to the equation $Z^2 = X^2 + 2Y^2$. Can be parametrically described in terms of three parameters:

$$\left\{\begin{array}{l} x = \delta \cdot \left|m^2 - 2n^2\right|, Y = 2\delta mn, Z = \delta(m^2 + 2n^2) \\ Where\ \delta,\ m,\ n\ are\ positive\ integers\ such \\ that\ gcd\ (m,\ n)\ =\ 1(i,e,m\ and\ n\ are\ relatively\ prime) \end{array}\right\} \quad \textbf{(14)}$$

## 7. Integral Volume and Proposition 1

Going back to **(1),** $\frac{a}{c} = \frac{b}{d} \Leftrightarrow \frac{a-c}{c} = \frac{b-d}{1}$; and so

$$\frac{d}{c} = \frac{b-d}{a-c} \quad \textbf{(15)}$$

From **(11)** and **(15)** we further obtain, we further obtain,

$$\frac{d}{c} = \frac{x}{y} \quad \textbf{(16)}$$

Now, we combine **(16), (11),** and **(8)** (the volume formula) to get,

$$V=\frac{H\cdot x}{3y}\cdot\left[\left(2y+\frac{\mathrm{d}y}{x}\right)^2+\left(2y+\frac{\mathrm{d}y}{x}\right)\frac{\mathrm{d}y}{x}+\left(\frac{\mathrm{d}y}{x}\right)^2\right];$$

$$V=\frac{H\cdot x}{3y}\cdot\frac{\left[y^2\cdot(2x+\mathrm{d})^2+y^2\cdot(2x+\mathrm{d})\cdot\mathrm{d}+\mathrm{d}^2\bullet y^2\right]}{x^2};$$

$$V=\frac{H\cdot y\cdot\left[(2x+d)^2+\mathrm{d}(2x+d)+d^2\right]}{3x}\text{ ; or equivalently}\qquad \textbf{(17)}$$

$$V=H\cdot y\cdot\left[2xd+\frac{3d^2+4x^2}{3x}\right];\text{ or equivalently}$$

$$V=2dyH+H\mathrm{d}\left(\frac{dy}{x}\right)+\frac{4Hyx}{3}$$

Recall from **(11)** that *1≤x≤y* and that x is a divisor of $y\cdot d$. Therefore, looking at the third version of the volume formula in **(17);** we see that the first and second terms are positive integers. The third term in **(17);** will be a positive integer precisely when *3* is a divisor of *Hyx*. And so, when this happens; the volume is an integer; and thus, the regular truncated pyramid with rectangular bases; is an integral one per **Definition 1**. We have the following:

***Proposition 1***

*Suppose that a, b, c, d, x, y, t, and H; are positive integers satisfying all the conditions in* ***(11).*** *Furthermore assume that the integer 3 is a divisor of Hyx.*

*Then the regular truncated pyramid with rectangular bases; and with base dimensions a, b (for the larger/bottom base) ; (for the smaller/top base) c, d; height H and the four lateral edges having length t; is an integral one (per* ***Definition 1****).*

## 8. Proposition 2 and Its Proof

In the end of the previous section, **Proposition 1** gives a general parametric description of the entire family of integral regular truncated pyramids with rectangular bases. In this section, we make use of **(13)** and **Proposition 1**; in order to obtain a 3-parameter subfamily of such pyramids. This is done in **Proposition 2** below.

***Proposition 2***

*Let l, m, and D be positive integer satisfying the following conditions:*

$$\left\{\begin{array}{l} M \text{ is a divisor of } l^2-1; \\ \text{That is, } l^2-1=m\cdot v, \\ \text{Where } v \text{ is a positive integer}; \\ 1\le m<l \\ l(m+v)\equiv 0(\bmod 3) \\ \end{array}\right\}$$

*Define the following positive integers d, x, y, a, c, b, d, H, and t:*

$$\begin{array}{l} d=\ 2mD,\ x=2m,\ y=2l,\ c=2lD, \\ a=c+2y=2l(D+2),\ b=d+2x=2m(D+2), \\ H=l^2+m^2-1=m(m+v),\ \text{and } t=l^2+m^2+1=m^2+mv+2 \end{array}$$

Then, the regular truncated pyramid with larger (rectangular) base lengths a and b; smaller (rectangular) base lengths c and d; height H, and with lateral edge length t; is a non-square integral (regular truncated) pyramid (with rectangular bases)

**Proof**

First note that the divisibility condition in **Proposition 1** is satisfied. Indeed, $yH=m(m+v)\cdot 2l=yH=2ml(m+v)$. And $3x=6m$.

Since $l\left(m+v\right)=0(\bmod 3)$; we have $l\left(m+v\right)=3i$, for some positive integer *i*. Thus,

$yH=2m\cdot 3i=6mi=i(6m)=2i\cdot\left(3x\right)$, which shows that the integer *3x* is a divisor of *yH.*

Next, we show that all the conditions in **(11)** are satisfied. Since $l>m$; $2l>2m$; $y=2l>2m=x$, so the condition $x\leq y$ is satisfied. The definitions of $a=c+2y$, $b=d+2x$; immediately imply that $a>c$ and $b>d$; by virtue of the fact that x and y are positive integers. Next compare a with b:

$a=2l(D+2)$ and $b=2m(D+2)$, which establishes $a>b$; in view $m<l$. Likewise

$c=2lD>d=2mD.$ This establishes that the regular truncated pyramid with rectangular bases is non-square. Finally the positive quadruple $(x,y,H,t)=(2m,\ 2l,\ m^2+l^2-1,\ m^2+l^2+1)$ satisfies the condition $x^2+y^z+H=t^2$ in **(11)**; by **(13)**. The proof is complete □.

Two remarks:

***Remark 1:*** If *m* is a divisor of $l-1$; then obviously $m$ is a divisor of $l^2-1=(l-1)(l+1)$; and clearly (since $m$ and $l$ are positive integers), $1\leq m\leq l-1<l$. Likewise, if $m$ is a divisor of $l+1$; then $m$ is a divisor of $l+1$; and so $m$ is a divisor of $l^2-1$ as well. Add to this the conditions $1\leq m\leq l+1$ and $l\geq 2$. Which clearly imply $1\leq m\leq l$; but then $m$ cannot equal $l$ (for $l$ cannot be a divisor of $l+1$; unless $l=1$ ). And so, $1\leq m<l$ is satisfied.

***Remark 2:*** The condition $l(m+v)=0(\bmod 3)$, says that 3 must divide at least one of $l, m+v$. If 3 is a divisor of $l$; then $l^2-1=mv$ implies that neither m nor v can be divisible by 3. And in that case, $mv\equiv l^2-1\equiv 0-1; mv\equiv 2(\bmod 3)$; which means that either $m\equiv 1$ and $v\equiv 2$; or alternatively $m\equiv 2$ and $v\equiv 1$(mod3). In either case, $m+v\equiv 0(\bmod 3)$. So, we see that when 3 divides $l$. Then the condition $l^2-1\equiv mv$; implies that 3 must also divide the sum $m+v$.

On the other hand, if $l$ is not divisible by 3. Then, $l\equiv 1$ or 2(mod3). And so $l^2\equiv 1$(mod3); $l^2-1\equiv 0(\bmod 3)$. Which implies that 3 divides the product $mv$. But since 3 also divides $l(m+v)$; and 3 does not divide $l$. It follows that 3 must divide $m+v$; and since it also divides $mv$. It follows that when $l$ is not a multiple of 3; both $m$ and $v$ are multiples of 3.

## 9. The Square Case and Proposition 3

When each of the two rectangular bases is a square; we have in **(11)**,

$$\left\{\begin{array}{l} a=b>c=d, y=x, \ and \ t^2=H^2+2y^2 \\ and \ a=2y+d \end{array}\right\} \quad \textbf{(18)}$$

Combining **(18)** with **(14)** produces,

$$\left\{\begin{array}{l} H=D\cdot\left|m^2-2n^2\right|, \ y=x=2Dmn, \ t=D(m^2+2n^2); \\ where \ m, \ n, \ D \ are \ positive \ integers; \ with \ m \ and \ n \\ being \ relatively \ prime. \end{array}\right\} \quad \textbf{(19)}$$

Combining **(19)** and **(18)** leads us to Proposition 3.

***Proposition 3***

*The set of integral regular truncated square pyramids with base lengths $a$ and $c$; $a>c$; height $H$, and lateral side length $t$; can be described in terms of four integer parameters as follows: $a=c+4Dmn$, $c=c$, $H=D\left|m^2-2n^2\right|$, and $t=D(m^2+2n^2)$. Where c, D, m, n are positive integers, with $m$ and $n$ being relatively prime.*

## 10. Closing remarks

The material on the subject matter of this work, seems pretty scant, as far as the literature is concerned. There is some material (a few pages long) on rational pyramids and trihedral angles, in L.E. Dickson's book, *History of The Theory of Numbers, Volume II*. See **Reference [2].**